\documentclass[11pt, leqno]{article}
\usepackage{amsmath, amsfonts, amsthm}
\begin{document}


\title{A symmetric diophantine equation\\ involving biquadrates}
\date{}
\author{Ajai Choudhry}

\maketitle

\theoremstyle{plain}

\newtheorem{lem}{Lemma}
\newtheorem{thm}{Theorem}

\begin{abstract} 
This paper is concerned with the diophantine equation $\sum_{i=1}^na_ix_i^4= \sum_{i=1}^na_iy_i^4$ where $n \geq 3$ and $a_i,\,i=1,\,2,\,\ldots,\,n$, are arbitrary integers. While a method of obtaining numerical solutions of such an equation has recently been given, it seems that an explicit parametric of this diophantine equation has not yet been published. We  obtain a multi-parameter  solution of this equation for arbitrary values of $a_i$ and for any positive integer $n \geq 3$, and deduce specific solutions when $n=3$ and $n=4$. The numerical solutions thus obtained are much smaller than the integer solutions of such equations obtained earlier.  
\end{abstract}

Keywords: biquadrates, fourth powers, quartic diophantine equation.

Mathematics Subject Classification 2010: 11D25
\bigskip
\bigskip

This paper is concerned with the diophantine equation
\begin{equation}
\sum_{i=1}^na_ix_i^4= \sum_{i=1}^na_iy_i^4, \label{quarteqnn}
\end{equation}
where $n$ is a positive integer  $\geq 3$ and $a_i,\,i=1,\,2,\,\ldots,\,n$, are arbitrary nonzero integers. We will be particularly interested in the following two special cases of Eq.~\eqref{quarteqnn}: 
\begin{align}
a_1x_1^4+a_2x_2^4+a_3x_3^4&=a_1y_1^4+a_2y_2^4+a_3y_3^4, \label{quarteqn3}\\
{\rm and} \quad a_1x_1^4+a_2x_2^4+a_3x_3^4+a_4x_4^4&=a_1y_1^4+a_2y_2^4+a_3y_3^4+a_4y_4^4. \label{quarteqn4}
\end{align}
 
Izadi and Baghalaghdam \cite{IB} have recently given a  method of obtaining numerical solutions of  the above three  equations by first relating their integer solutions to rational points of an elliptic curve. If the elliptic curve is of positive rank, infinitely many integer solutions of  these equations can be obtained.  It is, however, to be noted that for each of the three equations \eqref{quarteqnn}, \eqref{quarteqn3} and \eqref{quarteqn4}, if the elliptic curve concerned is of rank 0, the method given in \cite{IB} does not yield infinitely many  integer solutions of the equation under consideration.

It seems  that for arbitrary integer values of $a_i$, explicit parametric solutions of Eq.~\eqref{quarteqnn}, Eq.~\eqref{quarteqn3} and Eq.~\eqref{quarteqn4} have not been published till now. However, in the special case when the integers $a_i$ are all 1,  both Eq.~\eqref{quarteqn3} and Eq.~\eqref{quarteqn4} have been considered individually as well as part of a larger diophantine system, and several parametric solutions of  these diophantine equations are to be found scattered through the existing literature (see, for instance, \cite[pp.\ 305-306]{Ch1}, \cite[pp.\ 781-782]{Ch2}, \cite[pp.\ 653-657]{Di}).

We note that all the three equations \eqref{quarteqnn}, \eqref{quarteqn3} and \eqref{quarteqn4} are  homogeneous equations, and hence any  solution in rational numbers of these equations yields, on appropriate scaling, a primitive solution in integers. Therefore it suffices to find rational solutions of these equations. 

In this paper  we obtain a multi-parameter solution of Eq.~\eqref{quarteqnn} for any positive integer $n \geq 3$ and for arbitrary integer values of $a_i,\,i=1,\,2,\,\ldots,\,n$.  Parametric solutions of Eq.~\eqref{quarteqn3} and of Eq.~\eqref{quarteqn4} may be obtained by taking $n=3$ and $n=4$ respectively in the solution of Eq.~\eqref{quarteqnn}. The  numerical solutions  of  Eq.~\eqref{quarteqn3} and of Eq.~\eqref{quarteqn4} obtained from these parametric solutions are much smaller than the solutions of these equations obtained in \cite{IB} by using elliptic curves.

We note that the aforementioned multi-parameter solution of Eq.\eqref{quarteqnn} yields only trivial solutions of Eq.~\eqref{quarteqnn} when $n=2$ and the question of finding integer solutions of the quartic diophantine equation,
\begin{equation}
a_1x_1^4+a_2x_2^4=a_1y_1^4+a_2y_2^4, \label{quarteqn2}
\end{equation}
for arbitrary integer  values of $a_1$ and $a_2$  remains an open problem. 

We will now obtain a solution of Eq.~\eqref{quarteqnn} by following a  general method  described in \cite{Ch1} for solving symmetric diophantine systems. We have made minor modifications in the method so that the  solution can be expressed concisely. 

We note that if there is any solution of \eqref{quarteqnn} in which $x_i=\pm y_i$ for some $i$, we can cancel out the terms $a_ix_i^4$ and $a_iy_i^4$ on either side of \eqref{quarteqnn} and, in effect, we  get a solution of an equation of type \eqref{quarteqnn} with just $n-1$  biquadrates on either side. We will therefore  obtain only such solutions  of Eq.~\eqref{quarteqnn} in which $x_i \neq  \pm y_i, \; i=1,\,2,\,\ldots,\,n$.

To solve Eq.~\eqref{quarteqnn}, we write,
\begin{equation}
x_i=(f_i+g_i)u+r_iv,\;\;y_i=(f_i-g_i)u+r_iv,\;i=1,\,2,\,\ldots,\,n, \label{subs1}
\end{equation}
where $f_i,\,g_i,\,r_i,\,u$ and $v$ are arbitrary parameters such that $r_i \neq 0$. Since $x_i-y_i=2g_iu$, we impose the conditions $u \neq 0,\;g_i \neq 0,\;i=1,\,2,\,\ldots,\,n$, so that the solutions that we obtain satisfy the conditions $x_i \neq \pm y_i,\;i=1,\,2,\,\ldots,\,n$. 

On substituting the values of $x_i,\,y_i$ given by \eqref{subs1} in Eq.~\eqref{quarteqnn} and transposing all terms to the left-hand side, we get
\begin{multline}
\left(\sum_{i=1}^na_if_ig_i(f_i^2+g_i^2)\right)u^3+\left(\sum_{i=1}^na_ig_ir_i(3f_i^2+g_i^2) \right)u^2v\\
+3\left( \sum_{i=1}^na_if_ig_ir_i^2 \right)uv^2+\left(\sum_{i=1}^na_ig_ir_i^3 \right)v^3=0. \quad \quad \quad \quad \label{quarteqnnsubs1}
\end{multline}

We will now choose the parameters $f_i,\,g_i$ such that the coefficients of $uv^2$ and $v^3$ in Eq.~\eqref{quarteqnnsubs1} become 0. Equating to 0 the coefficient of $v^3$ in \eqref{quarteqnnsubs1}, we get a linear equation in $g_i,\;i=1,\,2,\,\ldots,\,n$, whose complete solution in rational numbers may be written as
\begin{equation}
g_i=(p_i-p_{i+1})/(a_ir_i^3),\;i=1,\,2,\,\ldots,\,n, \label{valg}
\end{equation}
where  $p_i,\;i=1,\,2,\,\ldots,\,n$, are  arbitrary rational parameters such that $p_i \neq p_{i+1}$ for  $i=1,\,2,\,\ldots,\,n$ and $p_{n+1}=p_1$. We note that there can only be $n-1$ independent linear parameters in the solution of the linear equation obtained by equating to 0 the coefficient of $v^3$ in \eqref{quarteqnnsubs1}, and hence one of the parameters $p_i,\; i=1,\,2,\,\ldots,\,n$ is actually superfluous but we prefer to write the solution in the symmetric manner  given by \eqref{valg}. The condition $p_i \neq p_{i+1}$ for each $i$ ensures that $g_i \neq 0,\;i=1,\,2,\,\ldots,\,n$.

With the values of $g_i$ now given by \eqref{valg},  on equating to 0 the coefficient of $uv^2$ in Eq.~\eqref{quarteqnnsubs1}, we get a linear equation in $f_i,\;i=1,\,2,\,\ldots,\,n$, whose complete solution in rational numbers may be written as
\begin{equation}
f_i=(q_i-q_{i+1})/(a_ig_ir_i^2),\;i=1,\,2,\,\ldots,\,n, \label{valf}
\end{equation}
where  $q_i,\;i=1,\,2,\,\ldots,\,n$ are  arbitrary rational parameters and $q_{n+1}=q_1$. Here again, one of the parameters $q_i,\; i=1,\,2,\,\ldots,\,n$ is  superfluous.

With the values of $f_i,\,g_i,\;i=1,\,2,\,\ldots,\,n$, being given by \eqref{valg} and \eqref{valf}, Eq.~\eqref{quarteqnnsubs1} reduces, on removing the nonzero factor $u^2$, to a linear equation in $u$ and $v$ whose solution is given by
\begin{equation}
u=\sum_{i=1}^na_ig_ir_i(3f_i^2+g_i^2),\quad v=-\sum_{i=1}^na_if_ig_i(f_i^2+g_i^2). \label{valuv}
\end{equation}

Thus a rational solution of Eq.~\eqref{quarteqnn} for arbitrary nonzero integer values of $a_i$ is given by \eqref{subs1} where $r_i,\; i=1,\,2,\,\ldots,\,n$, are arbitrary nonzero parameters, the values of $g_i,\; i=1,\,2,\,\ldots,\,n$, are given by \eqref{valg} in terms of arbitrary rational parameters  $p_i,\;i=1,\,2,\,\ldots,\,n$, such that $p_i \neq p_{i+1}$ for  $i=1,\,2,\,\ldots,\,n$ and $p_{n+1}=p_1$,  the values of  $f_i,\; i=1,\,2,\,\ldots,\,n$, are given by \eqref{valf} in terms of arbitrary rational parameters  $q_i,\;i=1,\,2,\,\ldots,\,n$ and  $q_{n+1}=q_1$, and the values of $u,\,v$ are given by \eqref{valuv}. 

A parametric solution of Eq.~\eqref{quarteqn3}, obtained  by taking $n=3,\,p_3=0,\,q_3=0$ and making suitable simplifications in the solution of Eq.~\eqref{quarteqnn}, is given by
\begin{equation}
x_i=(f_i+g_i)u+r_iv,\;\;y_i=(f_i-g_i)u+r_iv,\;i=1,\,2,\,3, \label{subs1eqn3}
\end{equation}
with the values of $f_i,\,g_i,\;i=1,\,2,\,3$ and $u,\,v$ being given by
\begin{equation}
\begin{aligned}
f_1 &= \lambda p_1p_2r_1(q_1+q_2), &g_1 &= -\mu a_2a_3r_2^3r_3^3(p_1+p_2),\\
f_2 &= \lambda p_1q_2r_2(p_1+p_2),  &g_2 &= \mu a_1a_3p_2r_1^3r_3^3, \\
f_3 &= \lambda p_2q_1r_3(p_1+p_2), &g_3 &= \mu a_1a_2p_1r_1^3r_2^3,\\
u&=\sum_{i=1}^3a_ig_ir_i(3f_i^2+g_i^2), & v&=-\sum_{i=1}^3a_if_ig_i(f_i^2+g_i^2), 
\end{aligned}
\end{equation}
where $p_1,\,p_2,\,q_1,\,q_2,\,r_1,\,r_2,\,r_3$ are arbitrary integer parameters while $\lambda$ and $\mu$ are arbitrary rational  parameters. 

As a numerical example, a solution of the equation,
\begin{equation}
x_1^4+x_2^4+61x_3^4=y_1^4+y_2^4+61y_3^4, \label{quarteqn3ex1}
\end{equation}
obtained by taking $p_1=61,\,p_2=-56,\,q_1=61,\,q_2=-21,\,r_1=1,\,r_2=2,\,r_3=-1,\, \lambda=1/4270,\,\mu=1/488$  is given by
\[6707^4+12802^4+61.(3237)^4=11227^4+6474^4+61.(4141)^4.\]

 This solution of Eq.~\eqref{quarteqn3ex1}  consists of  relatively small integers as compared to the smallest solution of Eq.~\eqref{quarteqn3ex1} given in \cite{IB} which involves large integers  having  24 digits.

 As in the case of Eq.~\eqref{quarteqn3}, we get a parametric solution of Eq.~\eqref{quarteqn4} by taking $n=4$ in the solution of Eq.~\eqref{quarteqnn} given  by \eqref{subs1}, \eqref{valg}, \eqref{valf} and \eqref{valuv}.  As a numerical example, a solution of the equation,
\begin{equation}
x_1^4+x_2^4+x_3^4+19x_4^4=y_1^4+y_2^4+y_3^4+19y_4^4, \label{quarteqn4ex1}
\end{equation}
obtained by taking $n=4,\, p_1=10,\,p_2=4,\,p_3=-4,\,p_4=-9,\,q_1=1,\,q_2=-2,\,q_3=6,\,q_4=1,\,r_1=1,\,r_2=2,\,r_3=1,\,r_4=-1$ in the parametric solution of Eq.~\eqref{quarteqnn}, is given by
\[576^4+220^4+527^4+19.(159)^4=600^4+416^4+453^4+19.(37)^4.\]
This solution is also much smaller than the smallest solution of Eq.~\eqref{quarteqn4ex1}, consisting of integers having 22 digits, given in \cite{IB}.

\begin{center}
\Large
Acknowledgment
\end{center}
 
I wish to  thank the Harish-Chandra Research Institute, Allahabad for providing me with all necessary facilities that have helped me to pursue my research work in mathematics.

\noindent Postal address: Ajai Choudhry, 13/4 A Clay Square, Lucknow - 226001, India\\
\noindent E-mail address: ajaic203@yahoo.com

\begin{thebibliography}{9}

\bibitem{Ch1} A. Choudhry, {\it Symmetric diophantine systems},   Acta
Arithmetica, \textbf{59} (1991), 291--307.

\bibitem{Ch2} A. Choudhry, {\it Equal Sums of Like Powers  and Equal Products of Integers},  Rocky Mountain Journal of Mathematics, \textbf{43} (2013),  763--792.

\bibitem{Di}L. E. Dickson, {\it History of theory of numbers, Vol. 2},
Chelsea Publishing Company, New York, 1992, reprint.

\bibitem{IB} F. Izadi and M. Baghalaghdam, {\it On the diophantine equation $\sum_{i=1}^nx_i^n= \sum_{j=1}^ny_j^n$}, available at  	arXiv:1701.02605

\end{thebibliography}
\end{document}